\documentclass[12pt]{article} 
\usepackage{amsfonts}
\usepackage{amsmath}
\usepackage{amssymb}
\usepackage{graphicx}
\usepackage{epsfig}
\date{\small\today}
\title{
An efficient solver for volumetric scattering \\based on fast spherical harmonics transforms}
\author{Youngae Han\\
\small Lorentz Solution, Inc.}
\newtheorem{theorem}{Theorem}[section]
\newtheorem{example}[theorem]{Example}

\begin{document}
\maketitle
\begin{abstract}
The Helmholtz equation arises in the study of electromagnetic radiation, optics, acoustics, etc.
In spherical coordinates, its general solution can be written as a spherical harmonic series which satisfies the radiation condition at infinity, ensuring that the wave is outgoing. 
The boundary condition at infinity is hard to enforce with a finite element method since a suitable approximation needs to be made within reasonable distance from scatterers.
 Luckily, the Helmholtz equation can be represented as a Lippmann-Schwinger integral equation which removes the necessity of the boundary approximations and its Green's function can be expanded as a spherical harmonic series which leads to our numerical scheme based on spherical harmonic polynomial transform. In this paper, we present an efficient solver for the Helmholtz equation which costs $O(N\log N)$ operations, where $N$ is the number of the discretization points. We use the fast spherical harmonic transforms which are originally developed in \cite{suda}. 
The convergence order of the method is tied to the global regularity of the solution. At the lower end, it is second order accurate for discontinuous material properties. The order increases with increasing regularity leading to spectral convergence for globally smooth solutions.
\end{abstract}

Lippmann-Schwinger integral equation, radiation condition, wave equation,\\
Spectral convergence.

\section{Introduction}
Computational electromagnetics and acoustics are fundamental to understanding of many practical systems like scattering, microwave circuits, radar, antennas, remote sensing, seismic exploration, ultrasound and tomography. Therefore, the demand for efficient numerical simulations of electromagnetic and acoustic fields is only increasing. The electromagnetic field is the solution of Maxwell's equations which are coupled with more than one unknown. These equations can be uncoupled by raising their order resulting in the wave equation. Thus, the Helmholtz equation reduced from the wave equation has been a classical problem (\cite{balanis}, \cite{jackson}) to solve and central focus of research for many decades and still drawing ongoing attentions. Using the separation of variables, the Helmholtz equation has been analytically solved for a simple geometry with homogeneous medium. In particular, for a spherical geometry, the solution is analytically given as the spherical harmonic series expansion in the angular direction and Bessel series expansion along the radial direction (\cite{balanis}) which inspires our numerical method based on spherical harmonic series expansion that is particularly efficient for spherical objects. There have been many efforts to develop fast polynomial transformations (\cite{healy}, \cite{fmm}, \cite{flt}, \cite{han_thesis}, \cite{suda}) to have equivalent speed advantages like FFTs. Here, we used the fast spherical harmonic transforms (\cite{han_thesis}) developed in \cite{suda} which is based on fast multipole method (\cite{fmm}) to have $O(N\log N)$ costs comparable to those based on FFTs while enjoying that the separation of variables that results from the addition theorem (\cite{colton}). This readily translates into quadratures that converge with higher order than those implicit in FFT-based schemes. \\
For the Helmholtz equation concerned in this paper is to find the scattered field generated due to the incident field with outgoing radiation condition. 
The most popular numerical methods are finite element methods (FEM) (\cite{goldstein}, \cite{harari}, \cite{ihlenburg}, \cite{lu}, \cite{thompson}) and integral equation methods (IEM) (\cite{anand_reitich}, \cite{bruno_hyde}, \cite{bruno_kunyansky}, \cite{han_thesis}, \cite{hyde_bruno}, \cite{rokhlin}, \cite{zhao}).  Although FEM can handle arbitrarily shaped obstacles with ease and the resulting matrix is sparse, it imposes a serious difficulty in solving scattering problems arising from the infinite size of domain.  For this reason, great effort has gone into the design of approximate local boundary conditions  and perfectly matched layers (\cite{berenger}, \cite{berenger_3d}, \cite{taeyoung}) that minimize spurious reflections which is by no means a trivial matter (\cite{harari}). 
Integral equation methods (IEM), on the other hand, implicitly account for radiation conditions through the use of \emph{outgoing} Green's functions. This very use of singular Green's functions, on the other hand, also translates into numerical challenges.  Moreover, these methods lead to a linear system involving a full matrix and thus they are not competitive unless a specialized strategy is used to accelerate matrix-vector products; examples of accelerated IEM include those based on FFTs (\cite{bruno_hyde}, \cite{hyde_bruno}, \cite{philips}) and those that use fast multipole expansions (\cite{fmm_lowfreq}, \cite{ling}, \cite{pan}).\\
 In this paper, inspired by the work in (\cite{bruno_sei}, \cite{bruno_hyde}, \cite{hyde_bruno}), we present a new accelerated IEM based on the addition theorem (\cite{colton}) and fast spherical harmonic transforms (\cite{han_thesis}, \cite{suda}). The convergence rate of our new algorithms is tied to the global regularity of fields. In particular, they converge with second order for the most singular case of discontinuous material properties and with increased rates for more regular arrangements; for smooth configurations the convergence is spectral.
This paper is organized as follows. The Lippmann-Schwinger integral equation is presented in Section 2. The numerical factorization of the integral equation based on spherical harmonic series expansion is explained in Section 3.  In Section 4, numerical implementations and their expected costs are derived, and in Section 5, numerical examples are given to confirm the predicted performance of the algorithms described in Section 4. Finally, in Section 6, the content of this paper is summarized.
\section{Lippmann-Schwinger integral equation}\label{sec2}
The wave equation states that
\begin{equation}\label{1}
 \frac{\partial^2 u}{\partial t^2}=\frac{c_0^2}{n^2}\nabla ^2 u 
\end{equation}
where $n(r)$ is the refractive index and $c_0$ is the propagation speed of the wave in air.
If we assume that the wave function is time-harmonic as  $$u(r,t)=u(r)e^{-iwt},$$ 
then a spatial solution $u(r)$ satisfies the Helmholtz equation
 \begin{equation}\label{2}
 \Delta {u}+k^2n^2{u} = 0 \mbox{ in } {\mathbb R ^3}.
\end{equation}
Here we will consider the scattering problem to determine the total field generated by a given incident field $u^i$.
Then the total field $ {u}$ which is the sum of  $u^i$ and the scattered field $u^s$ satisfies equation (\ref{2}) while $u^i$ is the solution of 
\begin{equation}\label{3}
 \Delta {u^{i}}+k^2{u^{i}} = 0 \mbox{ in } {\mathbb R ^3}.
\end{equation}
The \emph{Sommerfeld radiation condition} is given at infinity which guarantees that
the scattered field is outgoing,
\begin{equation}\label{4}
 \lim\limits_{r \rightarrow \infty} r\, \bigg (\frac{\partial u^{s} }{\partial r}-ik u^{s} \bigg)=0. 
\end{equation}
This boundary condition must be approximated at the reasonable distance from scatterers which is the main challenge to solve the Helmholtz equation.
But, this problem can be avoided if one appeals to the equivalent $\textit{Lippmann-Schwinger integral equation}$ which states 
\begin{equation}\label{5}
u(x)= u^{i}(x)-k^{2}\int\limits_{\Omega}^{}\! \Phi(x,y) u(y)m (y)\,dy,\ x \in \mathbb R^3
\end{equation}
where $m=1-n(x)^2$ and the inhomogeneity is of compact support $\Omega$. Therefore, we will develop an efficient solver based on the Lippmann-Schwinger integral equation.

\section{Spherical harmonic series expansion}\label{sec3}
To solve the integral equation (\ref{5}), we resort to the spherical harmonic series expansions of the Green's function $ \Phi(x,y)$. The addition theorem in \cite{colton} states that 
\begin{equation}\label{22}
\begin{split}
\Phi(x,y)&=\frac{1}{4\pi}\frac{e^{ik\mid x-y\mid}}{\mid x-y\mid}\\
&=ik\sum_{n=0}^{\infty}\sum_{m=-n}^{n}h_n^{(1)}(k\rho_>)
Y_n^m(\hat{\rho_>})j_{n}(k\rho_<)\stackrel{\underline{\hspace{0.5in}}}{Y_n^m(\hat{\rho_<})}
\end{split}
\end{equation}
where $$x=\widetilde{\rho}(\sin \widetilde{\theta}\cos \widetilde{\varphi},
\sin \widetilde{\theta} \sin \widetilde{\varphi}, \cos \widetilde{\theta} ),$$ $$y=\rho (\sin \theta \cos
\varphi,\sin \theta \sin \varphi,\cos \theta ),$$ 
$$\rho_<=\text{min}(\rho, \widetilde{\rho})  \text{ and } \rho_>=\text{max}(\rho, \widetilde{\rho}).$$ 
So, if we approximate $u$ and $u^i$ by a truncated spherical harmonic series as
 \begin{equation}\label{80}
  u^F(\rho, \theta,\phi)\}=\sum_{n=0}^{F}\sum_{m=-n}^{n} u^m_{n}(\rho)Y_n^m(\theta,\phi)
\end{equation}
\begin{equation}\label{79}
  u^{i,F}(\rho, \theta,\phi)\}=\sum_{n=0}^{F}\sum_{m=-n}^{n} u^{i,m}_{n}(\rho)Y_n^m(\theta,\phi)
\end{equation}
then from equation (\ref{5}) and the orthogonality properties of the spherical harmonics, $m(r)$ can be approximated without losing accuracy as
\begin{equation}\label{81}
  m^{2F}(\rho, \theta,\phi)\}=\sum_{n=0}^{2F}\sum_{m=-n}^{n} m^m_{n}(\rho)Y_n^m(\theta,\phi).
\end{equation}
Therefore, the formulation in (\ref{5}) becomes

\begin{equation}\label{78}
u_n^m(\rho)=u_n^{i,m}(\rho)
+iK_{n}^m(\rho),  
\end{equation}
where
\begin{equation}\label{77}
 K_{n}^m(\rho)=-k^3\int\limits_{0}^{\mathbf R}\!
h_{n}^{(1)}(k\rho_{>}) j_{n}(k\rho_{<}) I_{n}^m(\rho)\rho^2\,d\rho,
\end{equation}
and
\begin{equation}\label{76}
 I_{n}^m(\rho)
=\int\limits_{\phi=0}^{2\pi}\int\limits_{\theta=0}^{\pi}\!
u^{F}(\rho,\theta,\phi) m ^{2F}(\rho,\theta,\phi) 
\stackrel{\underline{\hspace{0.5in}}}{Y_{n}^m(\theta,\phi)}\sin\theta\,d\theta d\phi,
\end{equation}
$n=0,1,\dots,F$.
 To solve the formulation in (\ref{78}), we use the linear solver GMRES which requires the fast evaluation of angular integration of $I_{n}^m(\rho)$ in (\ref{76}) and the radial evaluation of $K_{n}^m(\rho)$ in (\ref{77}).

\section{Numerical implementation}\label{sec4}
\subsection{Angular integration}\label{sec6}
Due to the orthogonality of spherical harmonics, equation (\ref{76})
can be written as
\begin{equation}\label{75}
   u^F(\rho,\theta,\phi)m^{2F}(\rho,\theta,\phi)=\sum_{n=0}^{3F}\sum_{m=-n}^{n}
I_n^m(\rho)Y_n^m(\theta,\phi).
\end{equation}
Therefore, we define spherical harmonic transform as 
\begin{equation}\label{81}
   \{f(\theta_i,\phi_j)\}_{i,j=0}^{2F+1} \rightarrow \{\{c_n^m\}_{m=-n}^{n}\}_{n=0}^{F}
\end{equation}
where
$$
 f(\theta_i,\phi_j)=\sum_{n=0}^{F}\sum_{m=-n}^{n}c_n^mY_n^m(\theta_i,\phi_j)
$$
and its inverse, for appropriate choices of the angles $\{\theta_i,\phi_j\}$;  see \cite{suda}, \cite{han_thesis}.\\
In \cite{suda}, it is shown that if we define $S^{m}_n$ as 
$$
Y_n^m(t,\phi)=S^{|m|}_n(t)e^{im\phi}
$$
then
\begin{equation}\label{116}
\sum_{n=0}^{N-1}c_n^mS_{m+n}^m(x_j^N), \mbox{ j}\in \{0,1,2, ,N-1\}
\end{equation}
can be computed in $O(N\log N)$ operations for arbitrary $x_{j}^N$.\\
This work leads to fast spherical harmonic transform (FSHT) and its inverse (IFSHT) and these polynomial transforms cost $O(F^2\log F)$  for (\ref{81}).
 Thus, $I_n^m(\rho)$ of (\ref{76}) is computed as follows
\begin{equation}\label{50}
\begin{split}
\{\{I_{n}^m(\rho)\}_{m=-n}^{n}\}_{n=0}^{F}
&=FSHT_{3F}(IFSHT_{3F}(\{\{u_{n}^m(\rho)\}_{m=-n}^{n}\}_{n=0}^{F})\\ 
&\cdot IFSHT_{3F}(\{\{m_{n}^m(\rho)\}_{m=-n}^{n}\}_{n=0}^{2F})).
\end{split}
\end{equation}
Therefore, the angular integration costs $O(F^2\log F)$  for each $\rho$.

\subsection{Radial integration}\label{sec9}
The radial integral $K_{n}^m(\rho)$ in (\ref{77}) has a corner-type singularity at
$\rho=\widetilde{\rho}$  therefore, we write it as 
\begin{equation}\label{55}
\begin{split}
 \frac{-K_n^m(a)}{k^3}
 & =h_{n}^{(1)}(ka)\int\limits_{0}^{\min(a,\mathbf R)}\!j_{n}(k\rho)
   I_n^m(\rho)\rho^2\,d\rho
  +j_{n}(ka)\int\limits_{\min(a,\mathbf R)}^{\mathbf R}\!h_{n}^{(1)}(k\rho) 
   I_n^m(\rho)\rho^2\,d\rho \\
 &=i\bigg[ y_{n}(ka)\int\limits_{0}^{\min(a,\mathbf R)}\!j_{n}(k\rho)
   I_n^m(\rho)\rho^2\,d\rho
  +j_{n}(ka)\int\limits_{\min(a,\mathbf R)}^{\mathbf R}\!y_{n}(k\rho) 
   I_n^m(\rho)\rho^2\,d\rho \bigg] \\
 & +j_{n}(ka)\int\limits_{0}^{\mathbf R}\!j_{n}(k\rho)I_n^m(\rho)\rho^2\,d\rho. \\
\end{split} 
\end{equation}
Although it is numerically challenging to evaluate the Hankel function $h_n^{(1)}(k\rho)$, the product of the 
$j_{n}(k\rho)$  and  $h_n^{(1)}(k\rho)$ is bounded which leads us to define modified Bessel functions 
$\widetilde{j_n}(\rho)\mbox{ and }\widetilde{y_n}(\rho)$ in the form
\begin{equation}\label{35}
\begin{split}
&\widetilde{j_n}(\rho):=\frac{1\cdot3\cdot5\cdot\dots(2n+1)}{\rho^n}j_n(\rho)=
\bigg[1-\frac{\frac{1}{2}\rho^2}{1!(2n+3)}+\frac{({\frac{1}{2}\rho^2})^2}{2!(2n+3)(2n+5)}+\dots
\bigg]\\
&\widetilde{y_n}(\rho):=\frac{\rho^{n+1}}{-1\cdot1\cdot3\cdot5\cdot\dots(2n-1)}y_n(\rho)=
\bigg[1-\frac{\frac{1}{2}\rho^2}{1!(1-2n)}+\frac{({\frac{1}{2}\rho^2})^2}{2!(1-2n)(3-2n)}+\dots
\bigg]\\
\end{split}
\end{equation}
and with these modified Bessel functions we obtain 
\begin{equation}\label{83}
\begin{split}
K_n^m(a)
 &=\frac{i}{2n+1}\bigg[ \widetilde{y_n}(ka)\int\limits_{0}^{\min(a,\mathbf R)}\!
\big(\frac{\rho}{a}\big)^{n+1}\widetilde{j_n}(k\rho)
   I_n^m(\rho)k^2\rho\,d\rho \\
 &+\widetilde{j_n}(ka)\int\limits_{\min(a,\mathbf R)}^{\mathbf R}\!
\big(\frac{a}{\rho}\big)^{n}\widetilde{y_n}(k\rho) 
   I_n^m(\rho)k^2\rho\,d\rho \bigg] \\
 & +\widetilde{j_n}(ka)(ka)^n\int\limits_{0}^{\mathbf R}\!
\frac{(k\rho)^n(-2n-1)k^3}{1\cdot3^2\cdot5^2\dots (2n+1)^3}\widetilde{j_n}(k\rho)I_n^m(\rho)\rho^2\,d\rho. \\
\end{split} 
\end{equation}
For the radial integration, we divide the integration domain in a number $N_{i}$ of equi-length interpolation intervals $U_{j}=[u_j^0,u_j^1]$, $1\le j \le N_i$ on which we approximate for $\rho \in U_j,$ 
\begin{equation}\label{37}
I_{n}^{m}(\rho)\thickapprox\sum_{l=0}^{N_d-1}c_{l, m, n}^j T_{l}^{u_j^0,\mbox{}u_j^1}(\rho),
\end{equation}
where
$$
 T_{l}^{u_j^0,\mbox{}u_j^1}(\rho)=  T_{l}\bigg(\frac{\rho-(u_j^1+u_j^0)/2}{(u_j^1-u_j^0)/2}\bigg)
$$
are the Chebyshev polynomials in $U_{j}$. It costs $O(N_d (\log N_d)N_i)$ to compute the Chebyshev coefficients of $I_n^m(\rho)$ for fixed $n$ and $m$.
If we denote $\{\rho^j_k\}^{N_d}_{k=1}$ the Chebyshev points in $U_{j},$ the following equation holds
\begin{equation}\label{34}
\begin{split}
\int\limits_{0}^{\rho^j_{k+1}}\rho^{n+2} \widetilde{j_n}(k\rho)I_n^m(\rho) \,d\rho 
&=\int\limits_{0}^{\rho^j_{k}}\rho^{n+2} \widetilde{j_n}(k\rho)I_n^m(\rho) \,d\rho \\
& +\sum_{l=0}^{N_d-1}c_{l, m, n}^j  \int\limits_{\rho^j_{k}}^{\rho^j_{k+1}} \rho^{n+2} \widetilde{j_n}(k\rho) T_{l}^{u_j^0,\mbox{}u_j^1}(\rho) \,d\rho.\\
\end{split} 
\end{equation}
Therefore, if the moments $\int\limits_{\rho^j_{k}}^{\rho^j_{k+1}} \rho^{n+2} \widetilde{j_n}(k\rho) T_{l}^{u_j^0,\mbox{}u_j^1}(\rho)$ are computed and stored,
then it costs $O(N_d^2 N_i)$  to compute $\int\limits_{0}^{\rho^j_{k}}\rho^{n+2} \widetilde{j_n}(k\rho)I_n^m(\rho) \,d\rho$ for $1\le j \le N_i$ and  $1\le k \le N_d.$
Similarly, the computation of $\int\limits_{\rho^j_{k}}^{\mathbf R}\rho^{1-n} \widetilde{y_n}(k\rho)I_n^m(\rho) \,d\rho$ also costs $O(N_d^2 N_i)$ for $1\le j \le N_i$ and  $1\le k \le N_d.$ 
With these moments and precomputed coefficients $C_1(a)$, $C_2(a)$ and $C_3(a)$, $K_n^m(a)$ can be written as 
 \begin{equation}\label{80}
\begin{split}
K_n^m(a)
 &=C_1(a)\int\limits_{0}^{\min(a,\mathbf R)}\!
\rho^{n+2}\widetilde{j_n}(k\rho) I_n^m(\rho) \,d\rho \\
 &+  C_2(a)\int\limits_{\min(a,\mathbf R)}^{\mathbf R}\!
\rho^{1-n}\widetilde{y_n}(k\rho) 
   I_n^m(\rho) \,d\rho  \\
 & +C_3(a) \int\limits_{0}^{\mathbf R}\!
\rho^{n+2}\widetilde{j_n}(k\rho) I_n^m(\rho) \,d\rho. \\
\end{split} 
\end{equation}
Therefore the total cost for radial integration is
$$O\big[N_iN_d F^2 ( \log N_d+ N_d+1)\big] = O(N).$$
\section{Numerical examples}\label{sec5}
To show the predicted performance of the algorithm, four examples are presented below.
\begin{example}
\rm{Consider the scattering off a homogeneous sphere of radius 1 and the refractive index $n(x)=2$.
In this case, the problem is explicitly solvable, therefore comparison with the exact solution is possible.
A plane wave incident in the positive $z$-direction can be written as 
\begin{equation}\label{39}
e^{ik\vec{x}\cdot(0,0,1)}=e^{ik\rho \cos\theta}=\sum_{n=0}^{\infty}i^n(2n+1)j_{n}(k\rho)
 P_n(\cos\theta).
\end{equation}
If we differentiate equation (\ref{39}) $m_{inc}$ times with respect to $t=\cos\theta$,
we obtain the incident wave extended from (\ref{39}) as
\begin{equation}\label{86}
\begin{split}  
u^{i}&=\rho^{|m_{inc}|} e^{ik\rho \cos\theta} \sin^{|m_{inc}|} \theta e^{im_{inc}\phi}\\
&=\sum_{n=|m_{inc}|}^{\infty}
\frac{i^n(2n+1)}{(ik)^{|m_{inc}|}}j_{n}(k\rho)\frac{\sqrt{4\pi(n+|m_{inc}|)!}}{\sqrt{(2n+1)(n-|m_{inc}|)!}}
 Y_n^{m_{inc}}(\theta,\phi).\\
\end{split}
\end{equation}
For this incidence the exact solution \cite{balanis} is 
$$
u^{s}=
\begin{cases}
\sum_{n=|m_{inc}|}^{\infty}\{a_{n}j_{n}(2k\rho)-Q_n^{|m_{inc}|}(k,\rho)\}
Y_n^{m_{inc}}(\theta,\phi), 
&\rho \le1\\
\sum_{n=|m_{inc}|}^{\infty}b_{n}h_{n}(k\rho) Y_n^{m_{inc}}(\theta,\phi), &\rho \ge 1\\
\end{cases}
$$
where 
$$
  Q_n^{|m_{inc}|}(k,\rho)=\frac{i^n(2n+1)}{(ik)^{|m_{inc}|}}j_{n}(k\rho)\frac{\sqrt{4\pi(n+|m_{inc}|)!}}{\sqrt{(2n+1)(n-|m_{inc}|)!}}.
$$
By enforcing $C^1$ continuity of $u^s$ along the material discontinuity, $\{ a_n \}$ and $\{ b_n \}$ are determined.
The numerical error is computed as $E=\|u_{exact}-u_{approx}\|_{\infty}$ between the approximate and exact solutions for $m_{inc}=1$ and $k=5$ in
(\ref{86}) and different values of the interpolation orders $N_d$  in (\ref{37}) in Tables 1-3.}
\end{example}
\begin{figure}
\begin{center}
\fbox{
\includegraphics[width= 8cm]{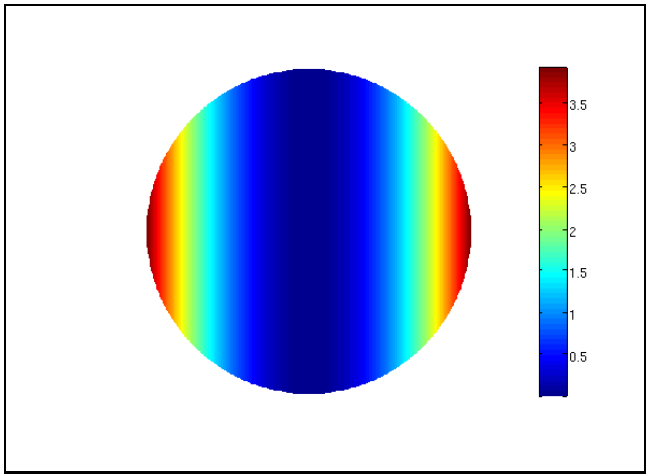}
}
\caption{The incident field intensity $|u^i|^2$ where $m_{inc}=1$ for examples 5.1-5.3.}
\label{sphere}
\end{center}
\end{figure}

\begin{figure}
\begin{center}
\fbox{
\includegraphics[width= 8cm]{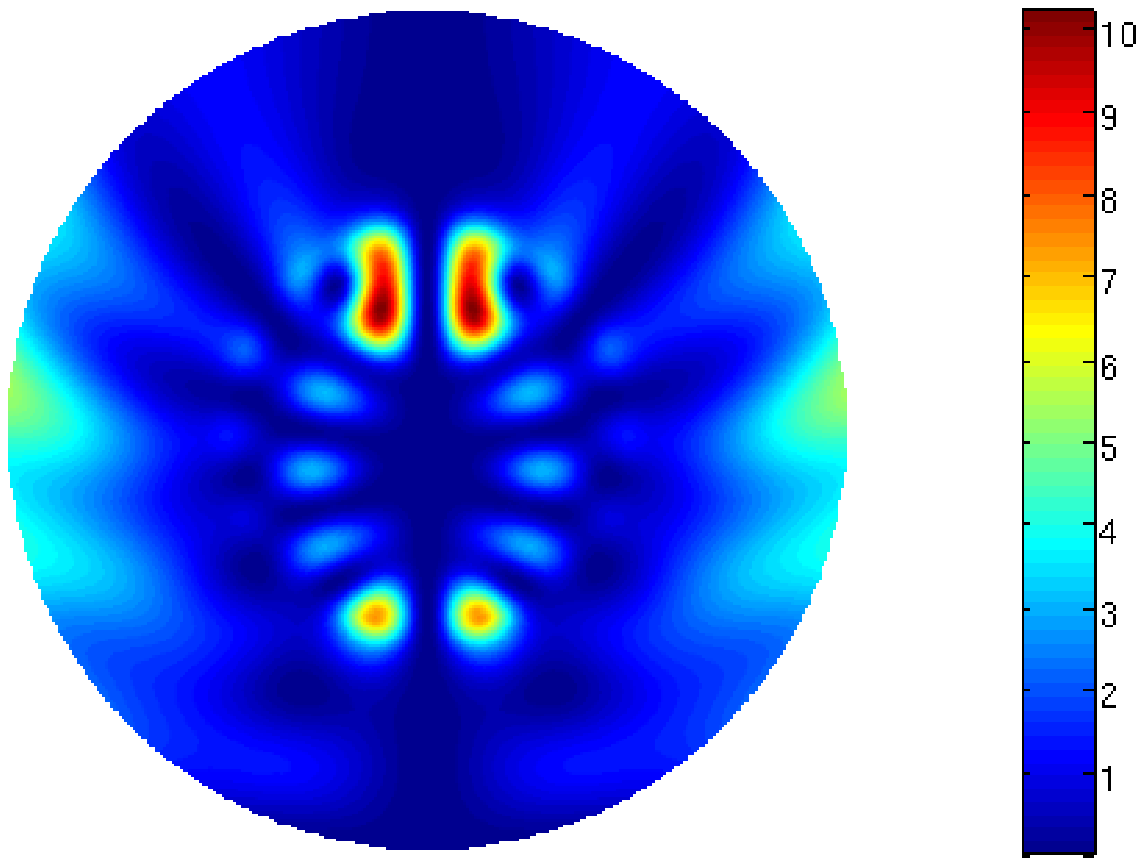}
}
\caption{The field intensity $|u|^2$ for example 5.1 where $m_{inc}=1$, $k =5$.}
\label{sphere}
\end{center}
\end{figure}

\begin{center}
\renewcommand{\arraystretch}{1.25}
\begin{tabular}{|r|p{1.6in}|p{1.5in} |r  |r|}
\hline
\textbf{$N_i$}  & \textbf{time per iteration}&\textbf{GMRES iteration}&\textbf{relative error}&\textbf{error ratio}\\ \hline
$2^3$ & 54 (sec)  &  34& 0.564482  & \\ \hline
$2^4$ & 108 (sec) &  35 & 0.20477  & 2.81207\\ \hline
$2^5$ & 217 (sec) &   36& 0.0532706  &  3.7178  \\ \hline
\end{tabular}
\end{center}
Table 1: Radial convergence for Example 5.1: the sphere centered at the origin.\\ 
$\mbox{  }\mbox{  }\mbox{  }\mbox{  }\mbox{  }\mbox{  }\mbox{  }\mbox{  }\mbox{  }\mbox{  }\mbox{  }\mbox{  }$
Parameters:  $m_{inc}=1,\ k=5,\ F=2^5-1,\ N_d=2,\ 0\le \rho \le 2,$\\ 
$\mbox{  }\mbox{  }\mbox{  }\mbox{  }\mbox{  }\mbox{  }\mbox{  }\mbox{  }\mbox{  }\mbox{  }\mbox{  }\mbox{  }$
GMRES tolerance = $1e-10.$\\

\begin{center}
\renewcommand{\arraystretch}{1.25}
\begin{tabular}{|r|p{1.6in}|p{1.5in} |r  |r|}
\hline
\textbf{$N_i$}  & \textbf{time per iteration}&\textbf{GMRES iteration}&\textbf{relative error}&\textbf{error ratio}\\ \hline
$2^3$ & 108 (sec) &  35& 0.00192818  & \\ \hline
$2^4$ & 217 (sec) &  36& 0.000193606& 10.0197\\ \hline
$2^5$ & 435 (sec) &  36& 1.35701e-05& 14.1395 \\ \hline
\end{tabular}
\end{center}
Table 2:  Radial convergence for Example 5.1 for $N_d=4$, GMRES tolerance = $1e-10,$\\
$\mbox{  }\mbox{  }\mbox{  }\mbox{  }\mbox{  }\mbox{  }\mbox{  }\mbox{  }\mbox{  }\mbox{  }\mbox{  }\mbox{  }$
Same parameters as in table 1.\\
\newpage
\begin{center}
\renewcommand{\arraystretch}{1.25}
\begin{tabular}{|r|p{1.6in}|p{1.5in} |r  |r|}
\hline
\textbf{$N_i$}  & \textbf{time per iteration}&\textbf{GMRES iteration}&\textbf{relative error}&\textbf{error ratio}\\ \hline
$2^3$ & 218 (sec)&  69 & 2.92336e-08 & \\ \hline
$2^4$ & 435 (sec)&  67 & 1.70495e-10 &  171.523 \\ \hline
$2^5$ & 1098 (sec)&  66 &  6.61144e-13  &   257.404 \\ \hline
\end{tabular}
\end{center}
Table 3:  Radial convergence for Example 5.1 for $N_d=8$, GMRES tolerance = $1e-15,$\\
$\mbox{  }\mbox{  }\mbox{  }\mbox{  }\mbox{  }\mbox{  }\mbox{  }\mbox{  }\mbox{  }\mbox{  }\mbox{  }\mbox{  }$
Same parameters as in table 1.\\

\begin{example}
\rm{To test the convergence when there is a material discontinuity for the angular integration, we consider the same sphere as in example 5.1 but now centered at (0, 0, $d$). For this case, we can obtain the incident wave as
\begin{equation}\label{87}
\begin{split}
u^{i}&=(x^2+y^2)^{\frac{|m_{inc}|}{2}}e^{ik(z-d)}e^{im_{inc}\phi}\\
&=e^{-ikd}\sum_{n=|m_{inc}|}^{\infty}Q_n^{|m_{inc}|}(k,\rho)
 Y_n^{m_{inc}}(\theta)\\
\end{split}
\end{equation}
and the exact spherical harmonic series expansion for $m(x)$ as follows
\begin{equation*}
\begin{split}
m(\rho,&\theta)\\
&=(1-n_0^2)\Big(\sum_{n=1}^{\infty}\frac{2n+1}{2}\sqrt{(1-\cos^2\theta_0)}\frac{(n-1)!}{(n+1)!}
P_n^1(\cos\theta_0)\sqrt{\frac{4\pi}{2n+1}}Y_n^0(\theta,\phi)\\
&+ \frac{1}{2}(1-\cos\theta_0)\sqrt{4\pi}Y_0^0(\theta,\phi)\Big)\\
\end{split}
\end{equation*}
where $0 \le(d-r) \le\rho \le(d+r)$, $n_0=2$ is the refractive index and 
$(\rho\sin\theta_0,\mbox{  }\rho\cos\theta_0)$ is a solution of 
$$
y^2+z^2=\rho^2,\mbox{  }y^2+(z-d)^2=r^2.
$$ 
The exact solution is obtained from the spherical harmonic transform of the exact values computed from the solution in example 5.1 with shifting.
In tables 4-6 we present radial and angular convergence studies for this case where the discontinuity does not lie on the grid. \\
}
\end{example}

\begin{center}
\renewcommand{\arraystretch}{1.25}
\begin{tabular}{|r|p{1.6in}|p{1.5in} |r  |r|}
\hline
\textbf{$N_i$}  & \textbf{time per iteration}&\textbf{GMRES iteration}&\textbf{relative error}& \textbf{error ratio}\\ \hline
$2^5$ &  1203 (sec) & 26  &0.299451  & \\ \hline
$2^6$ &  2408 (sec) & 27 & 0.0817522&3.64575 \\ \hline
$2^7$ &  4987 (sec) & 27 &  0.020442& 3.99036\\ \hline
\end{tabular}
\end{center}
Table 4:  Radial convergence for Example 5.2: the sphere centered
   at (0,0,2) with radius 1.\\ 
$\mbox{  }\mbox{  }\mbox{  }\mbox{  }\mbox{  }\mbox{  }\mbox{  }\mbox{  }\mbox{  }\mbox{  }\mbox{  }\mbox{  }$
Parameters:  $m_{inc}=1,\ k=5,\ F=2^7-1,\ N_d=2,\ 0\le \rho \le 4,$\\
$\mbox{  }\mbox{  }\mbox{  }\mbox{  }\mbox{  }\mbox{  }\mbox{  }\mbox{  }\mbox{  }\mbox{  }\mbox{  }\mbox{  }$
GMRES tolerance = $1e-5.$\\

\begin{center}
\renewcommand{\arraystretch}{1.25}
\begin{tabular}{|r|p{1.6in}|p{1.5in} |r  |r|}
\hline
\textbf{$N_i$}  & \textbf{time per iteration}&\textbf{GMRES iteration}&\textbf{relative error}& \textbf{error ratio}\\ \hline
$2^4$ & 601 (sec) & 6 &  2.88611e-05& \\ \hline
$2^5$ &  1203 (sec) & 6 & 8.11132e-06 & 3.51951 \\ \hline
$2^6$ &  2406 (sec) & 6 &  1.97643e-06& 4.0819\\ \hline
\end{tabular}
\end{center}
Table 5: Radial convergence for Example 5.2: the sphere centered
   at (0,0,2) with radius 1.\\ 
$\mbox{  }\mbox{  }\mbox{  }\mbox{  }\mbox{  }\mbox{  }\mbox{  }\mbox{  }\mbox{  }\mbox{  }\mbox{  }\mbox{  }$
Parameters:  $m_{inc}=3,\ k=1,\ F=2^7-1,\ N_d=2,\ 0\le \rho \le 4,$\\
$\mbox{  }\mbox{  }\mbox{  }\mbox{  }\mbox{  }\mbox{  }\mbox{  }\mbox{  }\mbox{  }\mbox{  }\mbox{  }\mbox{  }$
GMRES tolerance = $1e-10.$\\

\begin{center}
\renewcommand{\arraystretch}{1.25}
\begin{tabular}{|r|p{1.6in} |p{1.5in}| r    |r|}  
\hline
 \textbf{F}  & \textbf{time per iteration}& \textbf{GMRES iteration} &\textbf{relative error} & \textbf{error ratio}  \\ \hline
$2^4$-1 & 352 (sec)& 21 & 1.9425  &   \\ \hline
$2^5$-1 & 866 (sec)& 26 &   0.113651  &  17.0918 \\ \hline
$2^6$-1 & 2061 (sec)& 27 &    0.00157294& 72.2536  \\ \hline
\end{tabular}
\end{center}
Table 6: Angular convergence for Example 5.2: the sphere centered
   at (0,0,2) with radius 1.\\  
$\mbox{  }\mbox{  }\mbox{  }\mbox{  }\mbox{  }\mbox{  }\mbox{  }\mbox{  }\mbox{  }\mbox{  }\mbox{  }\mbox{  }$
Parameters:  $m_{inc}=1,\ k=5,\ N_d=2,\ N_i=2^7,\ 0\le \rho \le 4,$\\
$\mbox{  }\mbox{  }\mbox{  }\mbox{  }\mbox{  }\mbox{  }\mbox{  }\mbox{  }\mbox{  }\mbox{  }\mbox{  }\mbox{  }$
GMRES tolerance = $1e-5.$\\

\begin{example}
\rm{To test the convergence of a non-spherical object, we consider a square rotated by 45 degree and axisymmetric along z direction. 
 The refractive index $n(r)$ is 2 in $|x+y|\le 1$ shown in Figure 3.
The incident wave is the same as in example 5.1. We present radial and angular convergence studies in tables 7-9.
}
\end{example}

\begin{figure}
\begin{center}
\fbox{
\includegraphics[width= 8cm]{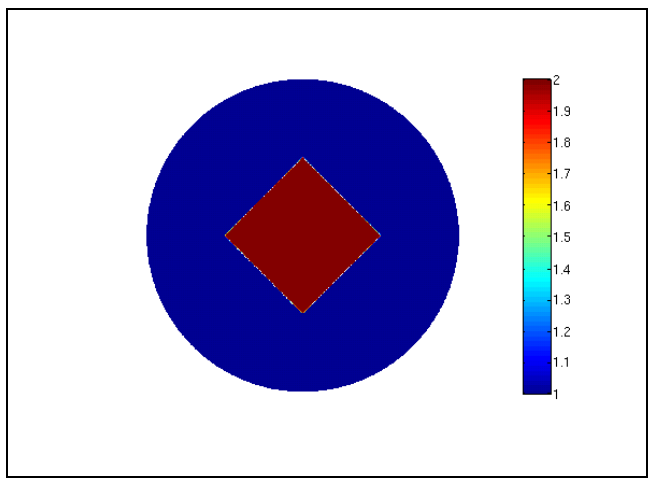}
}
\caption{The scatterer for example 5.3.}
\label{sphere}
\end{center}
\end{figure}

\begin{figure}
\begin{center}
\fbox{
\includegraphics[width= 8cm]{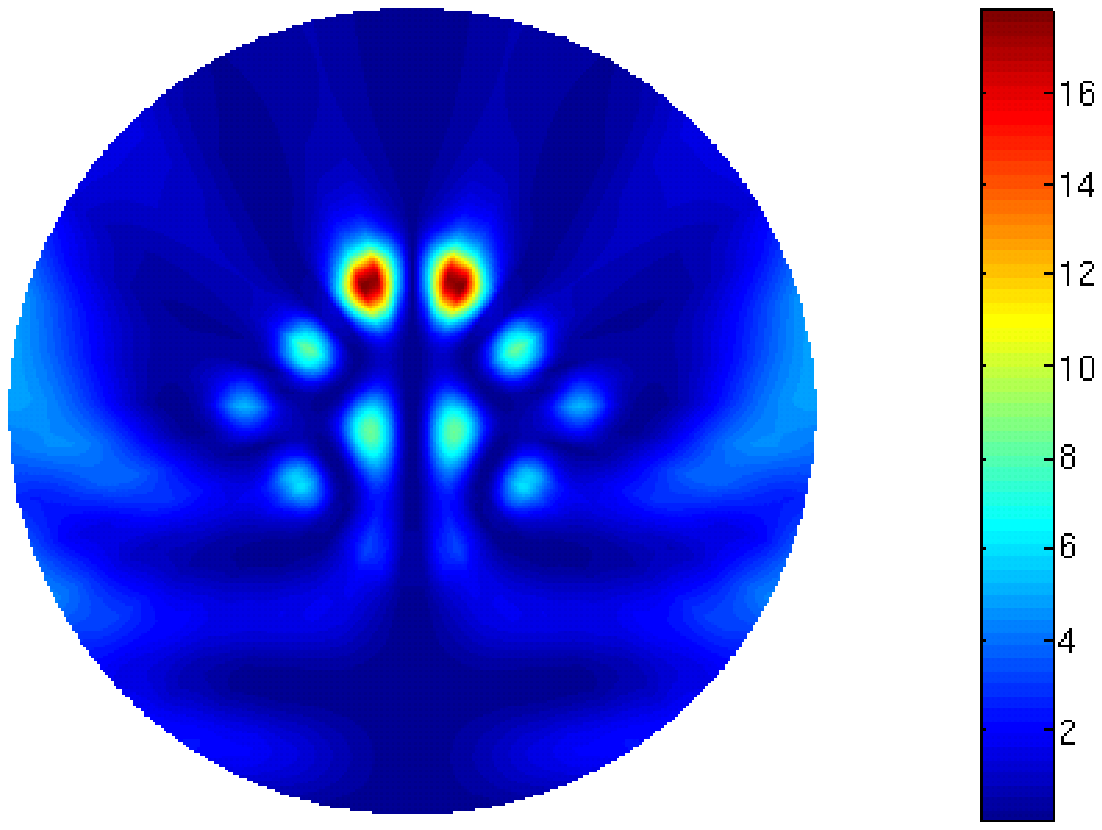}
}
\caption{The field intensity $|u|^2$ for example 5.3 where $m_{inc}=1$, $k =5$.}
\label{sphere}
\end{center}
\end{figure}

\begin{center}
\renewcommand{\arraystretch}{1.25}
\begin{tabular}{|r|p{1.6in}|p{1.5in} |r  |r|}
\hline
\textbf{$N_i$}  & \textbf{time per iteration}&\textbf{GMRES iteration}&\textbf{relative error}& \textbf{error ratio}\\ \hline
$2^4$ &  259 (sec) & 20 &    0.606676  & \\ \hline
$2^5$ &  559 (sec) & 20 &   0.0568001 & 10.6809\\ \hline
$2^6$ &  1081 (sec) & 20 &  0.0416413  & 1.36403 \\ \hline
$2^7$ &  2072 (sec) & 20 &   0.0126427 &  3.29371\\ \hline
\end{tabular}
\end{center}
Table 7:  Radial convergence for Example 5.3: the square rotated by 45 degree. \\ 
$\mbox{  }\mbox{  }\mbox{  }\mbox{  }\mbox{  }\mbox{  }\mbox{  }\mbox{  }\mbox{  }\mbox{  }\mbox{  }\mbox{  }$
Parameters:  $m_{inc}=1,\ k=5,\ F=2^6-1,\ N_d=2,\ 0\le \rho \le 2,$\\
$\mbox{  }\mbox{  }\mbox{  }\mbox{  }\mbox{  }\mbox{  }\mbox{  }\mbox{  }\mbox{  }\mbox{  }\mbox{  }\mbox{  }$
GMRES tolerance = $1e-5.$\\

\begin{center}
\renewcommand{\arraystretch}{1.25}
\begin{tabular}{|r|p{1.6in}|p{1.5in} |r  |r|}
\hline
\textbf{$N_i$}  & \textbf{time per iteration}&\textbf{GMRES iteration}&\textbf{relative error}& \textbf{error ratio}\\ \hline
$2^7$ &   2072 (sec) & 8 &   0.000389159   &  \\ \hline
$2^8$ &   4135 (sec) & 8 &   0.000123866 &3.14176\\ \hline
$2^9$ &  8729 (sec) & 8 &   1.79324e-05 & 6.90739\\ \hline
$2^{10}$ &  17823 (sec) & 8 &    3.94232e-06&  4.54869\\ \hline
\end{tabular}
\end{center}
Table 8:  Radial convergence for Example 5.3: the square rotated by 45 degree. \\ 
$\mbox{  }\mbox{  }\mbox{  }\mbox{  }\mbox{  }\mbox{  }\mbox{  }\mbox{  }\mbox{  }\mbox{  }\mbox{  }\mbox{  }$
Parameters:  $m_{inc}=1,\ k=1,\ F=2^6-1,\ N_d=2,\ 0\le \rho \le 2,$\\
$\mbox{  }\mbox{  }\mbox{  }\mbox{  }\mbox{  }\mbox{  }\mbox{  }\mbox{  }\mbox{  }\mbox{  }\mbox{  }\mbox{  }$
GMRES tolerance = $1e-10.$\\

\begin{center}
\renewcommand{\arraystretch}{1.25}
\begin{tabular}{|r|p{1.6in} |p{1.5in}| r    |r|}  
\hline
 \textbf{F}  & \textbf{time per iteration}& \textbf{GMRES iteration} &\textbf{relative error} & \textbf{error ratio}  \\ \hline
$2^4$-1 & 177 (sec)& 26&    0.0328253 &    \\ \hline
$2^5$-1 & 436 (sec)& 26 &   0.00439355 &    7.47126  \\ \hline
$2^6$-1 & 1036 (sec)& 26 &   0.000548512 &    8.00994 \\ \hline
\end{tabular}
\end{center}
Table 9: Angular convergence for Example 5.3: the square rotated by 45 degree. \\ 
$\mbox{  }\mbox{  }\mbox{  }\mbox{  }\mbox{  }\mbox{  }\mbox{  }\mbox{  }\mbox{  }\mbox{  }\mbox{  }\mbox{  }$
Parameters:  $m_{inc}=1,\ k=5,\ N_d=2,\ N_i=2^6,\ 0\le \rho \le 2,$\\
$\mbox{  }\mbox{  }\mbox{  }\mbox{  }\mbox{  }\mbox{  }\mbox{  }\mbox{  }\mbox{  }\mbox{  }\mbox{  }\mbox{  }$
GMRES tolerance = $1e-10.$\\

\begin{example}
\rm{The final example is to show the dependence of the scheme on the regularity of the scattering medium. We consider 
a refractive index given by
\begin{equation}\label{88}
n(\rho,\cos\theta)=
\begin{cases}
(1+\mid \cos\theta \mid ^{\beta}\sin^{|m_{ref}|}\theta e^{im_{ref}\phi})^{1/2}, &1\le \rho \le2\\
1, &0\le \rho <2 \mbox{ or } \rho >2.
\end{cases}
\end{equation}
Then, 
$$m(\rho,t)=
\begin{cases}
-|t|^{\beta}(1-t^2)^{\frac{|m_{ref}|}{2}}e^{im_{ref}\phi}, &1\le \rho \le2\\
0, &0\le \rho <1 \mbox{ or } \rho >2
\end{cases}
$$
and the exact spherical harmonic series expansion for $m(\rho, t)$ is given by $\sum_{l=0}^{\infty}m_{2l}(\rho)Y_{|m_{ref}|+2l}^{m_{ref}}(t)$, 
where
\begin{equation}\label{89}
\begin{split}
m_{2l}(\rho)&=-\int_{\phi=0}^{2\pi}\int_{\theta=0}^{\pi}\!
Y_{|m_{ref}|+2l}^{m_{ref}}(\theta,\phi)\cos^{\beta}\theta
\sin^{|m_{ref}|}\theta e^{im_{ref}\phi}\cos\theta\,d\theta d\phi\\
&=-\sqrt{\frac{(2n+1)(n-|m_{ref}|)!}{(n+|m_{ref}|)!}}\frac{\pi}{2^{\beta+|m_{ref}|}}
\frac{\Gamma(1+\beta)}{\Gamma(1+\frac{\beta}{2})}
\frac{1}{\Gamma(\frac{3}{2}+\frac{\beta}{2})}\\
&\cdot \frac{\prod_{s=0}^{l-1}(\frac{1}{2}\beta-s)}
{\prod_{s=0}^{|m_{ref}|+l-1}(\frac{1}{2}\beta+\frac{3}{2}+s)} \frac{(n+|m_{ref}|)!}{(n-|m_{ref}|)!}.\\
\end{split}
\end{equation}
 In tables 10-12, we present the order of convergence for $\beta=$0.4, 1.4, 2.4 to show the correlation between smoothness 
of the refractive index $n(x)$ and the order of convergence.
} 

\begin{center}
\renewcommand{\arraystretch}{1.25}
\begin{tabular}{|r|p{1.6in} |p{1.5in}| r    |r|}  
\hline
 \textbf{F}  & \textbf{time per iteration}& \textbf{GMRES iteration} &\textbf{relative error} & $\log_2$(\textbf{error ratio})  \\ \hline
$2^3$-1 & 83 (sec)& 6& 0.00186872 &    \\ \hline
$2^4$-1 & 419 (sec)& 6& 7.42083e-06 &   7.97625\\ \hline
$2^5$-1 & 2048 (sec)& 6& 1.19402e-06 &   2.63575\\ \hline
\end{tabular}
\end{center}
Table 10: Angular convergence for Example 5.4 for
$\beta=0.4 \ (m\in C^{0.4},\ u\in C^{2.4}).$\\
$\mbox{  }\mbox{  }\mbox{  }\mbox{  }\mbox{  }\mbox{  }\mbox{  }\mbox{  }\mbox{  }\mbox{  }\mbox{  }\mbox{  }$
Parameters:  $m_{inc}=3,\ m_{ref}=1,\ k=0.5,\ N_d=8,\ N_i=4,\ 0\le \rho \le 4,$\\
$\mbox{  }\mbox{  }\mbox{  }\mbox{  }\mbox{  }\mbox{  }\mbox{  }\mbox{  }\mbox{  }\mbox{  }\mbox{  }\mbox{  }$
GMRES tolerance = $1e-10.$\\

\begin{center}
\renewcommand{\arraystretch}{1.25}
\begin{tabular}{|r|p{1.6in} |p{1.5in}| r    |r|}  
\hline
 \textbf{F}  &  \textbf{GMRES iteration} &\textbf{relative error} & $\log_2$(\textbf{error ratio})  \\ \hline
$2^3$-1 &  5 &  0.0018687 &    \\ \hline
$2^4$-1 &  5 &  6.99546e-07 &   11.3833 \\ \hline
$2^5$-1 &  5 &  5.82286e-08  &    3.58662\\ \hline
\end{tabular}
\end{center}
Table 11: Angular convergence for Example 5.4 for
$\beta=1.4 \ (m\in C^{1.4},\ u\in C^{3.4}).$\\
$\mbox{  }\mbox{  }\mbox{  }\mbox{  }\mbox{  }\mbox{  }\mbox{  }\mbox{  }\mbox{  }\mbox{  }\mbox{  }\mbox{  }$
Same parameters as in table 10.\\

\begin{center}
\renewcommand{\arraystretch}{1.25}
\begin{tabular}{|r|p{1.6in} |p{1.5in}| r    |r|}  
\hline
 \textbf{F}& \textbf{GMRES iteration} &\textbf{relative error} & $\log_2$(\textbf{error ratio})  \\ \hline
$2^3$-1 & 5 &  0.00186869  &    \\ \hline
$2^4$-1 & 5 &  6.6957e-08 &  14.7684 \\ \hline
$2^5$-1 & 5 &  2.74242e-09 &   4.60971 \\ \hline
\end{tabular}
\end{center}
Table 12: Angular convergence for Example 5.4 for
$\beta=2.4 \ (m\in C^{2.4},\ u\in C^{4.4}).$\\
$\mbox{  }\mbox{  }\mbox{  }\mbox{  }\mbox{  }\mbox{  }\mbox{  }\mbox{  }\mbox{  }\mbox{  }\mbox{  }\mbox{  }$
Same parameters as in table 10.\\

\begin{figure}
\begin{center}
\fbox{
\includegraphics[width= 8cm]{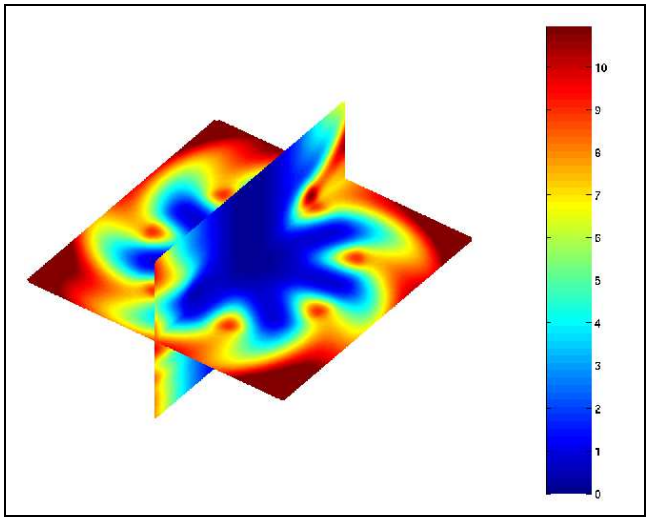}
}
\caption{ The field intensity $|u|^2$  for example 5.4 where $m_{inc}=1$, $m_{ref}=7$, $\beta=0.2$,   $k =5$. }
\label{sphere}
\end{center}
\end{figure}

\end{example}

\section{Summary}
In this paper, an efficient solver for scattering by penetrable three-dimensional structures is presented.
The solution is obtained by the iterative evaluation of Lippmann-Schwinger integral equation and 
its efficiency comes from the use of the addition theorem and fast spherical harmonics transforms. 
The scheme allows for such evaluations in $O(N\log N)$ operations, where $N$ is the number of the discretization points.
The convergence order of the method, on the other hand, is tied to the global
regularity of the solution. At the lower end, it is second order accurate for discontinuous material properties. The order increases with increasing regularity of the refractive index leading to spectral convergence for globally smooth solutions.

\section*{Acknowledgments}
Y. Han would like to thank Dr. E. M. Hyde for a helpful discussion. This work was in part supported by AFOSR Contract F49620-02-1-0052.

\end{document}